\documentclass[a4paper]{article}
\usepackage[margin=3cm]{geometry}
\usepackage{graphicx,xcolor} 
\usepackage{amsmath,amsthm,amsfonts,amssymb}
\usepackage{enumitem}
\usepackage{mathrsfs}

\newcommand{\curlyd}{\mathtt{d}}

\usepackage[colorlinks=true,
linkcolor=blue,
citecolor=green!50!black,  
urlcolor=teal]{hyperref}

\usepackage{biblatex}
\newtheorem{thm}{Theorem}
\newtheorem{lem}{Lemma}
\newtheorem{rmk}{Remark}
\newtheorem{defn}{Definition}

\title{On dual-projectively equivalent connections associated to second order superintegrable systems}
\author{\textsc{Andreas Vollmer}\thanks{{University of Hamburg, Department of Mathematics, 
			Bundesstra{\ss}e 55, 20146 Hamburg, Germany}, \texttt{andreas.vollmer@uni-hamburg.de} } }
\date{\small 5 May 2025}

\begin{document}
	
	\maketitle
	
	\begin{abstract}
		Pre-geodesics of an affine connection are the curves that are geodesics after a reparametrization (the analogous concept in K\"ahler geometry is known as $J$-planar curves).
		Similarly, dual-geodesics on a Riemannian manifold are curves along which the $1$-forms associated to the velocity are preserved after a reparametrization.
		
		Superintegrable systems are Hamiltonian systems with a large number of independent constants of the motion. They are said to be second order if the constants of the motion can be chosen to be quadratic polynomials in the momenta. Famous examples include the Kepler-Coulomb system and the isotropic harmonic oscillator.
		
		We show that certain torsion-free affine connections which are naturally associated to certain second order superintegrable systems share the same dual-geodesics.    
	\end{abstract}
	
	\section{Introduction}
	
	We consider geometric structures (more precisely, certain affine connections) that naturally appear in the context of second order (maximally) superintegrable Hamiltonian systems. Such systems include famous models from mathematical physics, such as the Kepler-Coulomb system, the isotropic harmonic oscillator or the Smorodinski-Winternitz system. We obtain that these structures are dual-projectively equivalent, a concept which has been introduced in the context of statistical manifolds, Weylian structures and affine hypersurfaces.
	
	Let $(M,g)$ be a Riemannian (smooth) manifold.
	Assume that, for $\varepsilon>0$, $\gamma:(-\varepsilon,\varepsilon)\to M$ is a (smooth) curve on $M$ with tangent (velocity) vector field $\dot\gamma$. We denote the $1$-form associated to $\dot\gamma$ (by virtue of $g$) by $\dot\gamma^\flat$.
	Here, $\flat:\mathfrak X(M)\to\Omega^1(M)$ denotes the usual musical isomorphism induced by $g$. Similarly, we denote by $\sharp:\Omega^1(M)\to\mathfrak X(M)$ the musical isomorphism induced by $g^{-1}$, when the underlying metric is clear.
	\begin{defn}[\cite{Ivanov1995}]
		A curve $\gamma$ on $M$ is called \emph{dual-geodesic} for an affine connection $\nabla$ if
		\[
		\nabla_{\dot\gamma}\dot\gamma^\flat = q(t)\,\dot\gamma^\flat,
		\]
		where $q:(-\varepsilon,\varepsilon)\to\mathbb R$.
		In particular, we say that $\gamma$ is an \emph{affinely parametrized} dual-geodesic for $\nabla$, if
		\[
		\nabla_{\dot\gamma}\dot\gamma^\flat = 0.
		\]
		If $\nabla$ is the Levi-Civita connection of the (Riemannian) metric $g$, then we also say that a curve is dual-geodesic for $g$, if it is dual-geodesic for $\nabla$.
	\end{defn}
	
	It is well-known that for every dual-geodesic curve, there exists an affine parametrization, see \cite[Prop.~2.1]{Ivanov1995}. In this reference, dual-geodesics are introduced as a tool for the study of semi-conjugate connections and affine hypersurface immersions, and we refer the interested reader there for more detailed information on this perspective.
	Here, we mention only the following fact, which we need later: \emph{Let $p\in M$ and $w\in T_pM$. Then there exists a (unique up to reparametrization and possibly with smaller $\varepsilon$) dual-geodesic curve $\gamma:(-\varepsilon,\varepsilon)\to M$, $\gamma(0)=p$ with $\dot\gamma(0)=w$}, see \cite[Prop.~2.2]{Ivanov1995} .
	
	\begin{defn}[\cite{Ivanov1995}]
		Two connections are called dual-projectively equivalent, if they share the same dual-geodesic curves.
	\end{defn}
	
	Dual-geodesics and dual-projectively equivalent connections have been discussed, for instance, in \cite{Ivanov1995,Matsuzoe2010}, where they have been related to  affine hypersurfaces, statistical manifolds and Weylian structures.
	The purpose of this paper is to demonstrate that dual-projectively equivalent connections naturally arise in the context of second order superintegrable Hamiltonian systems. Let $(M,g)$ be a simply connected (connected) Riemannian manifold and denote its Levi-Civita connection by~$\nabla$. Then $T^*M$ carries a natural symplectic structure $\omega$ induced by the tautological $1$-form.
	We consider a natural Hamiltonian $H:T^*M\to\mathbb R$,
	\[
	H(x,p)=g^{-1}_x(p,p)+V(x)\,,
	\]
	where $(x,p)$ are canonical Darboux coordinates on $T^*M$.
	For a function $f:T^*M\to\mathbb R$, we denote by $X_f$ the Hamiltonian vector field with respect to the natural symplectic structure, i.e.
	$
	\iota_{X_f}\omega=df.
	$
	
	\begin{defn}
		A (maximally) \emph{superintegrable system} is given by a Hamiltonian~$H$ together with $2n-2$ functions $F^{(m)}:T^*M\to\mathbb R$,
		such that $(H,F^{(1)}, \dots,$ $F^{(2n-2)})$ are functionally independent,
		and such that
		$
		X_H(F^{(m)})=0
		$
		for all $1\leq m\leq 2n-2$.
		We say that a superintegrable system is \emph{second order} if the functions $F^{(m)}$ are quadratic polynomials in the momenta, i.e.
		\[
		F^{(m)}(x,p)=K_{(m)}^{ij}(x)p_ip_j+W^{(m)}(x).
		\]
	\end{defn}
	For the integrals of motion in a second order superintegrable system, it is easy to check that (omitting the subscript $(m)$ for brevity) the tensor field $\sum g_{ai}g_{bi}K^{ab}\,dx^i\otimes dx^j$ is Killing, i.e.\ satisfies $\nabla_XK(X,X)=0$ for all $X\in\mathfrak X(M)$.
	We write $\mathcal K$ for the $\mathbb R$-linear space of Killing tensors associated to a second order superintegrable system, meaning that there is a function $W$ on $M$ such that $F=K(p^\sharp,p^\sharp)+W$ is an integral of the motion for $H$, i.e.\ $X_H(F)=0$.
	
	\begin{defn}
		We say that a second order superintegrable system is \emph{irreducible}, if the linear space generated by the endomorphisms $K$, $K^\flat\in\mathcal K$, form an irreducible set, i.e.\ do not share a common eigenspace.
		For the sake of brevity, an irreducible second order superintegrable system will simply be referred to as an \emph{irreducible system} in the following.
	\end{defn}
	
	It was proven in~\cite{KSV2023} that, for an irreducible system, there exists a tensor field $\hat T\in\Gamma(\mathrm{Sym}^2_\circ\otimes TM)$, trace-free in its covariant indices, such that
	\begin{equation}\label{eq:raw.wilczynski}
		\nabla^2V = \hat T(dV)+\frac1n\,g\,\Delta V,
	\end{equation}
	where $\hat T$ depends on $\mathcal K$ only, and where $\Delta$ denotes the Laplace-Beltrami operator. In general, the tensor $\hat T$ is not unique, but here we confine ourselves to systems for which $\hat T$ is unique. Specifically, we consider \emph{non-degenerate} second order superintegrable systems. These are irreducible systems with a $(n+2)$-parameter family of potentials (see Section \ref{sec:irreducible} for a precise definition). The main results are Theorems~\ref{thm:main.1} and~\ref{thm:main.2} in Section~\ref{sec:proof}, which show that three affine connections, which are naturally defined for non-degenerate systems, are dual-projectively equivalent:
	\begin{enumerate}[label=$\langle$\Alph*$\rangle$,noitemsep]
		\item\label{item:ind} the induced connection $\nabla^g\pm\hat T$ (``induced connections''),
		\item\label{item:stat} the corresponding connection that endows the space with the information-geometric structure of a \emph{statistical manifold},
		\item\label{item:semi} the connections that naturally arise when one restricts to an $(n+1)$-dimensional subspace of potentials (to be explained later).
	\end{enumerate}
	These connections can also be found in \cite{KSV2023,KSV2024,CV2025,BLMS24,NV}, for example.
	Before we prove these dual-projective equivalences, we review some facts about irreducible second order superintegrable Hamiltonian systems.

	\section{Irreducible second order superintegrable systems}\label{sec:irreducible}
	
	Two specific kinds of irreducible systems are going to play a crucial role in the following, namely non-degenerate and (generalised) semi-degenerate systems. These are introduced in the following two subsections. The terminology goes back to the foundational work by Kalnins and coworkers, cf.\ \cite{KKM18,KKM-I,KKM-II,KKM-III,KKM-IV,KKM-V} and the references therein. For semi-degenerate systems, we also mention \cite{Escobar-Ruiz&Miller}.

	\subsection{Non-degenerate systems}
	
	\emph{Non-degenerate systems} are quadruples $(M,g,\mathcal K,\mathcal V)$ such that $(M,g)$ is as before, $\mathcal K$ a linear space of Killing tensors (of dimension $2n-1$ or larger, with $g\in\mathcal K$) and a linear subspace $\mathcal V\subset\mathcal C^\infty(M)$ (of dimension $n+2$), such that the space of endomorphisms associated to $\mathcal K$ is irreducible, and 
	$ d(K(dV))=0 $
	for all $K\in\mathrm{End}(TM)$ such that $K^\flat\in\mathcal K$.
	Such a system satisfies~\eqref{eq:raw.wilczynski}, which then implies the (closed) prolongation system ($\nabla$ denotes the Levi-Civita connection of $g$, and $\Delta$ its Laplace-Beltrami operator)
	\begin{align*}
		\nabla^2V &= \hat T(dV)+\frac1n\,\Delta V\,g
		\\
		\nabla\Delta V &= \hat q(dV)+(\mathrm{tr}(\hat T)-q)\Delta V
	\end{align*}
	where $q(X,Y)=g(\hat q(X),Y)$ and
	$
	\hat q(X)=\mathrm{tr}_g(\nabla_{\cdot}\hat T(\cdot,X))+\mathscr T(X)-\widehat{\mathrm{Ric}}^g(X)
	$
	with $g(\widehat{\mathrm{Ric}^g}(X),Y)=\mathrm{Ric}^g(X,Y)$.
	Also, we introduce $\mathscr T\in\mathrm{End}(TM)$ via
	$\mathscr T(X)=\mathrm{tr}_g(\Theta(X,\cdot,\cdot))$, 
	where $\Theta:\mathfrak X(M)^3\to\mathfrak X(M)$,
	\[
	\Theta(X,Y,Z)(\alpha)=\hat T(X,Y)(\hat T(Z)(\alpha)),
	\]
	for $X,Y,Z\in\mathfrak X(M)$, $\alpha\in\Omega^1(M)$, where $\hat T(Z)(\alpha)$ is the $1$-form $\hat T(\cdot,Z)(\alpha)$.
	For a non-degenerate system, we define the induced connections by
	\begin{equation}\label{eq:nabla.T}
		\nabla^{\pm\hat T}:=\nabla\mp\hat T,
	\end{equation}
	which is torsion-free and Ricci-symmetric, see \cite{BLMS24}. For simplicity, we abbreviate $\nabla^{\hat T}=\nabla^{+\hat T}$.
	Following \cite{KSV2023}, we furthermore introduce the totally symmetric and tracefree tensor field $S\in\Gamma(\mathrm{Sym}^3_\circ(T^*M))$ and the $1$-form $t\in\Omega^1(M)$ by setting $S=\mathring T$ and $t=\frac{n}{(n-1)(n+2)}\mathrm{tr}(\hat T)$, such that
	\[
	T(X,Y,Z)=S(X,Y,Z)+t(X)g(Y,Z)+t(Y)g(X,Z)+t(Z)g(X,Y)\,,
	\]
	where $T:=T^\flat$ and where $X,Y,Z\in\mathfrak X(M)$.
	Next, for dimension $n\geq3$, we let
	\begin{equation}
		\mathcal{Z}=\mathscr S-(n-2)(S(t)+t\otimes t)-\mathrm{Ric}^g,
	\end{equation}
	where $\mathscr S(X,Y)=\mathrm{tr}_g(~\hat S(X,\cdot)(S(\cdot,Y))~)$.
	It is shown in \cite{KSV2023} that, if $n\geq3$ and if the underlying manifold is of constant sectional curvature, then
	\begin{equation}\label{eq:zeta}
		\mathring{\mathcal{Z}}=\mathring\nabla^2\zeta
	\end{equation}
	for a function $\zeta\in\mathcal C^\infty(M)$.
	We can hence introduce the totally symmetric tensor field
	\[
	\digamma=T+\frac{n+2}{n}\,g\otimes t+\frac1{2(n-2)}\,\Pi_{\mathrm{Sym}}g\otimes d\zeta\,.
	\]
	which is then also a Codazzi tensor, c.f.~\cite{KSV2023}.
	Note that the definition of $\digamma$ relies on the assumption of having a space of constant sectional curvature. For later use, we also introduce
	$\hat\digamma=\digamma g^{-1}$ and $\nabla^{\pm\hat\digamma}=\nabla\mp\hat\digamma$ (and $\nabla^{\hat\digamma}=\nabla^{+\hat\digamma}$).
	
	Relaxing the curvature assumptions again, we introduce, for a non-degenerate system in dimension $n\geq2$, the totally symmetric tensor field
	\[
	B=T+\frac{n+2}{n}\,g\otimes t\,,
	\]
	as well as the connections
	\[
	\nabla^{\pm\hat B}:=\nabla\mp{\hat B}\qquad (\nabla^{\hat B}=\nabla^{+\hat B}).
	\]
	where $\hat B=Bg^{-1}$, c.f.\ \cite{KSV2023}.
	%
		We remark that for so-called \emph{abundant} systems, the connections $\nabla^{\pm\hat\digamma}$ and $\nabla^{\pm\hat B}$, respectively, coincide up to a suitable \emph{gauge choice} of $\zeta$.
		An abundant system is a non-degenerate system with $\frac12n(n+1)$ linearly independent, compatible Killing tensor fields.
		Note the non-trivial freedom for choosing the function~$\zeta$,  satisfying~\eqref{eq:zeta}. This gauge freedom is thoroughly discussed in~\cite{KSV2023}.
		If $n\geq3$ and $g$ has constant sectional curvature, \cite{KSV2023} shows that one can choose $\zeta=0$ without changing the data of $S$ and $t$, i.e.\ without modifying the structure tensor~$\hat T$.

	\subsection{Semi-degenerate systems}
	
	\emph{Generalized semi-degenerate systems}, or $(n+1)$-parameter systems, are quadruples $(M,g,\mathcal K,\mathcal V)$ such that $(M,g)$ is as before, $\mathcal K$ a linear space of Killing tensors (of dimension $2n-1$ or larger, with $g\in\mathcal K$) and a linear subspace $\mathcal V\subset\mathcal C^\infty(M)$ (of dimension $n+1$), such that the space of endomorphisms associated to $\mathcal K$ is irreducible, and
	$ d(K(dV))=0. $
	Moreover, we require that in addition to~\eqref{eq:raw.wilczynski}, an equation of the form
	\begin{equation}\label{eq:semi-degeneracy}
		\Delta V=\hat s(dV)
	\end{equation}
	holds, for some $\hat s\in\mathfrak X(M)$ that is determined by $\mathcal K$, and where $\Delta$ is the Laplace-Beltrami operator of the Levi-Civita connection $\nabla$ of $g$.
	For the generalized semi-degenerate system subject to \eqref{eq:semi-degeneracy} we therefore have
	\begin{equation*}
		\nabla^2V = \hat D(dV)\,,
	\end{equation*}
	where we introduce the tensor field $\hat D=\hat T+\frac1n\,g\otimes\hat s\in\Gamma(\mathrm{Sym}^2(T^*M)\otimes TM)$.
	Note that $\hat D$ depends on the space $\mathcal K$ only. We also introduce $D=\hat D^\flat\in\Gamma(\mathrm{Sym}^2(T^*M)\otimes T^*M)$ for later use.
	
	We say that a generalized semi-degenerate system, is \emph{weak}, if there is $\mathcal V'\supset\mathcal V$ such that $(M,g,\mathcal K,\mathcal V')$ is non-degenerate.
	Otherwise, we call it a \emph{strong} semi-degenerate system.
	For a (weak or strong) semi-degenerate system, we define the induced connection by
	\[
	\nabla^{\pm\hat D}:=\nabla\mp\hat D
	\]
	(again abbreviating $\nabla^{\hat D}=\nabla^{+\hat D}$).
	It is shown in \cite{BLMS24} that $\nabla^{\hat D}$ is torsion-free, Ricci-symmetric and projectively flat (the reference only discusses the case of strong semi-degeneracy, but it is easy to extend this result to generalized semi-degenerate systems).
	We also introduce the tensor field, c.f.~\cite{NV},
	\begin{multline*}
		N(X,Y,Z) := \frac{1}{3}\left( 2\,D(X,Y,Z) -D(X,Z,Y) -D(Y,Z,X)\right)
		\\
		+ \frac{1}{3(n-1)} \ \left(
		2g(X,Y) \curlyd(Z)
		- g(X,Z) \curlyd(Y)
		- g(Y,Z) \curlyd(X)
		\right)
	\end{multline*}
	where
	$\curlyd = (n+2)t - s$.
	It is shown in \cite{NV} that $N=0$ characterizes precisely the situation of a generalized semi-degenerate system that is weak (i.e.\ it is strong if $N$ does not vanish).
	This means, in the case $N=0$, that
	\[
	\hat T = \hat D-\frac1n\,g\otimes\hat s
	\]
	satisfies the conditions of a non-degenerate structure tensor.
	For later use, and to keep the notation clean, we introduce the $1$-form $s\in\Omega^1(M)$, $s=\hat s^\flat$.

	\section{Proof of the main results}\label{sec:proof}
	
	In this section, we show the dual-geodesic equivalence of the connections \ref{item:ind}--\ref{item:semi}. All of these connections are torsion-free. Indeed, denoting the Levi-Civita connection of $g$ by $\nabla$, the torsion-freeness of $\nabla\pm\hat T$ follows immediately from the symmetries of $\hat T$, cf.~\cite{BLMS24}. The torsion-freeness of $\nabla\pm\hat B$ follows immediately from the total symmetry of $B^\flat$. In the semi-degenerate case, the torsion-freeness of the connections $\nabla\pm\hat D$ follows similarly.
	
	Before we proceed to the actual proof, we review some results from the literature that are going to be useful later.
	
	\begin{lem}[Prop.~2.3 of \cite{Ivanov1995}]\label{lem:formula}
	Let $(M,g)$ be a pseudo-Riemannian manifold.
	Then two torsion-free affine connections $\nabla,\nabla'$ are dual-projectively equivalent if and only if there is a $1$-form $\alpha\in\Omega^1(M)$ such that
	\begin{equation}\label{eq:dual-proj.equiv}
		\nabla'_XY=\nabla_XY+\alpha^\sharp\, g(X,Y)
	\end{equation}
	for any vector fields $X,Y\in\mathfrak X(M)$.
	\end{lem}
	\noindent Torsion-freeness is a necessary requirement for \eqref{eq:dual-proj.equiv}, and in the presence of torsion counterexamples can easily be found.
	
	For later use, we also introduce the concept of \emph{semi-compatibility} for pairs $(\nabla',h)$ consisting of an affine connection $\nabla'$ and a metric $h$.
	
	\begin{defn}[\cite{Ivanov1995}]
		The pair $(\nabla',h)$ is said to be \emph{semi-compatible} (via $\alpha$), if there exists a $1$-form~$\alpha$ such that
		\[
			\nabla'_Xh(Y,Z)-\nabla'_Yh(X,Z)=\alpha(Y)h(X,Z)-\alpha(X)h(Y,Z)
		\]
		for any $X,Y,Z\in\mathfrak X(M)$.
		The pair $(\nabla',h)$ is called \emph{compatible}, if it is semi-compatible via $\alpha=0$.
	\end{defn}
	
	We begin our investigation with the dual-projective equivalence of the connections~\ref{item:ind} and~\ref{item:stat}.	
	To this end, consider a non-degenerate system on $(M,g)$ with structure tensor $\hat T$ as before.
	Observe that the induced connection $\nabla^{\hat T}$ and the structural connection $\nabla^{\hat B}$ are dual-projectively equivalent.
	
	\begin{thm}\label{thm:main.1}
		\begin{enumerate}[label=(\roman*),wide,noitemsep]
			\item The connections~$\nabla^{\hat T}$ and~$\nabla^{\hat B}$ share the same dual-geodesics.
			\item For a non-degenerate system with induced connection $\nabla^{\hat T}$, there is a unique dual-projectively equivalent connection $\nabla^\star$ such that $(\nabla^\star,g)$ is compatible. In fact, $\nabla^\star=\nabla^{\hat B}$.
		\end{enumerate}
	\end{thm}
	\noindent The analogous statements hold, if we replace $\nabla^{\hat T}$ and $\nabla^{\hat B}$ by $\nabla^{-\hat T}$ and $\nabla^{-\hat B}$, respectively.
	We comment that the following proof also shows that $T=0$, if $\nabla^{\hat T}=\nabla^\star$. This latter condition holds for the so-called non-degenerate harmonic oscillator system \cite{KSV2023}.	
	\begin{proof}
		We denote the Levi-Civita connection of $g$ by $\nabla^g$.
		We have
		$\nabla^{\hat T}-\nabla^{\hat B}
		= \nabla^g-\hat T-\nabla^g+\hat B
		= \hat B-\hat T.$
		Using the musical isomorphisms, we then compute
		\begin{equation*}
			(B-T)(X,Y,Z)
			= \frac{n+2}{n}\,t(Z)g(X,Y)
		\end{equation*}
		and conclude
		$
		\nabla^{\hat T}-\nabla^{\hat B} = \frac{n+2}{n}\,g\otimes t^\sharp\,.
		$
		This proves the first claim.
		We next consider the connections that are dual-projectively equivalent to $\nabla^{\hat T}$. They are of the form, $\beta\in\Omega^1(M)$,
		\[
		\nabla^\star_XY=\nabla^{\hat T}_XY+\beta^\sharp\,g(X,Y).
		\]
		A short computation shows that
		\[
		\nabla^\star_Xg(Y,Z)-\nabla'_Yg(X,Z)
		= \alpha(Y)g(X,Z)-\alpha(X)g(Y,Z)
		\]
		with the $1$-form
		$
		\alpha=\frac{n+2}{n}t-\beta.
		$
		The connection $\nabla^\star$ therefore is compatible with $g$ if and only if
		$
		\beta = \frac{n+2}{n}t.
		$
		We conclude
		\[
		\nabla^\star_XY=\nabla^{\hat T}_XY+\frac{n+2}{n}t^\sharp\,g(X,Y)=\nabla^{\hat B}_XY.
		\]
	\end{proof}
	
	\begin{rmk}
		We remark that an analogous computation shows
		$
		\nabla^{\hat\digamma}_Xg(Y,Z)-\nabla^{\hat\digamma}_Yg(X,Z)
		= 0\,,
		$
		alongside
		$
		\nabla^{\hat B}_Xg(Y,Z)-\nabla^{\hat B}_Yg(X,Z)
		= 0\,.
		$
		However, the connections $\nabla^{\hat\digamma}$ and $\nabla^{\hat B}$ are, in general, different, as
		\[
		g(\nabla^{\hat\digamma}-\nabla^{\hat B})
		=\frac1{2(n-2)}\,\Pi_{\mathrm{Sym}}g\otimes d\zeta.
		\]
		We infer that the connections $\nabla^\digamma$ and $\nabla^{\hat B}$ coincide precisely if $d\zeta=0$.
		Note that the vanishing of $d\zeta$ implies $\mathring{\mathcal{Z}}=0$.
	\end{rmk}

	We now turn our attention to the dual-projective equivalence of the connections \ref{item:ind} and~\ref{item:semi}, i.e.\ we now consider systems with (n+1)-parameter potential.	Again, we focus on $\nabla^{\hat D}=\nabla^{+\hat D}$ for conciseness, as the discussion for $\nabla^{-\hat D}$ is analogous.
	We introduce the connection
	\[
	\nabla^{\dagger}=\nabla^{\hat D}-\frac1n\,s^\sharp\,g
	\]
	which is clearly dual-projectively equivalent to $\nabla^{\hat D}$.
	We characterize weak semi-degeneracy via~$\nabla^{\hat D}$.
	
	\begin{thm}\label{thm:main.2}
		\begin{enumerate}[label=(\roman*),wide,noitemsep]
			\item Consider a weak semi-degenerate system with induced connection~$\nabla^{\hat D}$.
			Assume that the induced connection of the naturally associated non-degenerate system is $\nabla^{\hat T}$.
			Then $\nabla^{\hat D}$ and $\nabla^{\hat T}$ share the same dual-geodesics.
			\item Consider a (generalized) semi-degenerate system with induced connection $\nabla^{\hat D}$ and semi-degeneracy $1$-form $s$.
			Then $(\nabla^{\hat D},g)$ are semi-compatible via
			\[ \beta = \frac1n\,\left( s-(n+2)\,t\right)\,, \]
			if and only if the system is a weak semi-degenerate system.
		\end{enumerate}
	\end{thm}
	\noindent The analogous statements hold, if we replace $\nabla^{\hat T}$ and $\nabla^{\hat D}$ by $\nabla^{-\hat T}$ and $\nabla^{-\hat D}$, respectively.
	\begin{proof}
		We have $\nabla^{\hat T}=\nabla^\dagger$, and hence $\nabla^{\hat T}$ and $\nabla^{\hat D}$ are dual-projectively equivalent, completing the first part of the theorem.
		For the second part, we first compute
		\begin{align}
			\nabla^{\hat D}_Xg(Y,Z)-\nabla^{\hat D}_Yg(X,Z)
			&= N(Y,Z,X)-N(X,Z,Y)
			\nonumber \\
			&\qquad + \frac1n\,\left( s(X)-(n+2)\,t(X) \right)\,g(Y,Z)
			\nonumber \\
			&\qquad - \frac1n\,\left( s(Y)-(n+2)\,t(Y) \right)\,g(X,Z)
			\nonumber \\
			&\stackrel{!}{=} \beta(X)g(Y,Z)-\beta(Y)g(X,Z)
			\label{eq:beta.condition}
		\end{align}
		where the exclamation point 
		indicates the requirement 
		that $(\nabla^{\hat D},g)$ be semi-compatible via~$s$.
		
		\noindent Part ``$\Rightarrow$'': Inserting the formula for $\beta$ into~\eqref{eq:beta.condition}, we obtain the condition
		$$
		N(Y,Z,X)=N(X,Z,Y).
		$$
		It follows that $N=0$ and, invoking~\cite{NV}, we thus obtain the claim.
		
		\noindent Part ``$\Leftarrow$'': If the system is weakly semi-degenerate, then $N=0$,  due to~\cite{NV}. We immediately find that the condition at the exclamation point holds, if $\beta$ is as claimed.
	\end{proof}

	\section{Conclusion}
	We have seen here that certain affine connections that naturally appear in the theory of irreducible superintegrable systems are dual-projectively equivalent.
	In particular, the theorems stated in this paper imply, that extendability (weak semi-degeneracy) for a $(n+1)$-parameter system is linked to the semi-compatibility (with the metric $g$) of its induced connection $\nabla^{\hat D}$. Weak semi-degeneracy in turn implies that there is a naturally associated non-degenerate system whose induced connection $\nabla^{\hat T}$ is dual-projectively equivalent to $\nabla^{\hat D}$. In this case there is also a connection $\nabla^{\hat B}$ that is compatible with $g$ and dual-projectively equivalent to $\nabla^{\hat D}$. 
	The observed occurrence of dual-projective geometry is natural and linked to the underlying Weylian structure. The underlying Weylian structure was discussed in~\cite{Vollmer}. Note that, by a direct computation,
	\[
	\nabla^{\hat T}_Xg(Y,Z)-\frac{n+2}{n}\,t(X)g(Y,Z) \in\Gamma(\mathrm{Sym}^3(T^*M))\,
	\]
	is totally symmetric. According to~\cite{Matsuzoe2010}, it was shown in~\cite{Matsuzoe2007} that this implies the existence of a Weylian connection.

	\section*{Acknowledgements}
	
	This research was funded by the German Research Foundation (Deutsche For\-schungs\-gemein\-schaft) through the research grant \#540196982.

	\AtNextBibliography{\small}
	\printbibliography

\end{document}